\documentclass[11pt]{article}
\usepackage{amsfonts}
\usepackage{color,xcolor}
\usepackage{amsmath}
\usepackage{latexsym,amsfonts,amssymb}
\usepackage{mathrsfs}

\newcommand{\vsss}{\vspace{.05in}}

\begin{document}

\title{ Asymptotic behavior of solutions of a free boundary problem modeling tumor spheroid with Gibbs--Thomson relation }

\author{Junde Wu$^{\dag}$\ \ \ and\ \ \ Fujun Zhou$^{\ddag}$
\\[0.2cm]
  {\small $^{\dag}$ Department of Mathematics, Soochow University,
  Suzhou, }\\ [-0.1cm]
  {\small Jiangsu 215006, PR China. E-mail:\,wujund@suda.edu.cn}
  \\[0.1cm]
  {\small $^{\ddag}$ Department of Mathematics, South China University of Technology,
  }
  \\  [0.1cm]
  {\small Guangzhou, Guangdong 510640, PR China. E-mail:\,fujunht@scut.edu.cn}
  \\ [-0.1cm]
}

\date{}

\medskip

\maketitle

\begin{abstract}
  In this paper we study a free boundary problem modeling the growth of solid tumor
  spheroid. It consists of two elliptic equations describing nutrient diffusion and pressure
  distribution within tumor, respectively. The new feature is that nutrient concentration
  on the boundary is less than external supply due to a Gibbs--Thomson relation and the problem has two
  radial stationary solutions, which differs from widely studied tumor spheroid model
  with surface tension effect. We first establish local well-posedness  by using a
  functional approach based on Fourier multiplier method and analytic semigroup theory.
  Then we investigate stability of each radial stationary solution. By employing a generalized
  principle of linearized stability, we prove that the radial stationary solution with
  a smaller radius is always unstable, and there exists a positive threshold value
  $\gamma_*$ of cell-to-cell adhesiveness $\gamma$, such that the radial stationary
  solution with a larger radius is
  asymptotically stable for $\gamma>\gamma_*$, and unstable for $0<\gamma<\gamma_*$.

\medskip


\medskip


\medskip

\end{abstract}

\section{Introduction}
\setcounter{equation}{0}
\hskip 1em

  In the last several decades, numerous tumor models have been established to explore the mechanism of tumor
  growth [\ref{low-10}, \ref{roose-07}]. Since the tumor region is evolving with time, Greenspan [\ref{gre-72}] first
  pointed out that it is natural to model tumor growth in form of free boundary problems of partial
  differential equations. Byrne and Chaplain developed Greenspan's idea
  and proposed the following model [\ref{byr-cha-96}, \ref{byr-cha-97}]:
\begin{equation}\label{1.1}
\left\{
\begin{array}{ll}
   \Delta\sigma=\sigma \qquad \qquad \qquad \qquad &
   \hbox{in  $\Omega(t) \times (0,+\infty) $},
   \\ [0.2 cm]
    \Delta p=-\mu(\sigma-\tilde{\sigma})\quad  \qquad  &\hbox{in  $\Omega(t)
    \times (0,+\infty)$},
\\ [0.2 cm]
  \sigma=\bar{\sigma}(1-\gamma\kappa) \qquad\qquad & \hbox{on
  $\partial\Omega(t) \times (0,+\infty)$},
\\ [0.2 cm]
   p=\bar{p} \qquad \qquad \qquad \qquad & \hbox{on $\partial\Omega(t)
   \times (0,+\infty)$},
\\ [0.2 cm]
     V=-\partial_{\bf n} p \qquad \qquad \qquad & \hbox{on $\partial\Omega(t)
     \times (0,+\infty)$},
\\ [0.2 cm]
   \Omega(0)=\Omega_0,
\end{array}
\right.
\end{equation}
  where $\sigma=\sigma(x,t)$ and $p=p(x,t)$ denote nutrient concentration  and internal pressure within
  the domain $\Omega(t)\subset \Bbb R^3$ occupied by tumor at time $t>0$, respectively. It is a free boundary
  problem, in which $\sigma(x,t)$, $p(x,t)$ and $\Omega(t)$ are unknowns and have to be determined together.
  $\kappa$ is the mean curvature, $V$ is the outward normal velocity and ${\bf n}$ denotes the unit outward normal of tumor boundary
  $\partial\Omega(t)$, respectively.  $\bar{\sigma}$, $\bar p$,
  $\tilde{\sigma}$, $\mu$ and $\gamma$ are positive constants, where $\bar{\sigma}$ represents constant
  external nutrient supply, $\bar p$ denotes constant external pressure, $\tilde{\sigma}$ is a critical
  nutrient concentration for apoptosis, $\mu$ is the proliferation rate of tumor cells and $\gamma$ is
  cell-to-cell adhesiveness on tumor boundary. $\Omega_0$ is the region occupied by tumor initially.

  The above model (\ref{1.1}), which is usually called {\em tumor spheroid model with Gibbs--Thomson relation},
  was proposed based on the following considerations [\ref{byr-cha-96}, \ref{byr-cha-97}]. Firstly, it is
  hypothesized that energy is expended in maintaining the tumor's compactness by cell-to-cell adhesion on the
  boundary and the nutrient acts as a source of such energy and satisfies Gibbs--Thomson relation. This relation
  means that nutrient concentration on the tumor boundary is less than external nutrient concentration and
  the jump is proportional to local mean curvature. Secondly, the pressure on the boundary is assumed to
  be equal to external pressure, reflecting that it is continuous across the tumor boundary. For a special
  simplification of this model by replacing the first equation of (\ref{1.1}) with $\Delta\sigma=C$ for a constant $C$,
  radial solution was studied and linear stability of radial stationary solution under asymmetric perturbations
  was obtained in [\ref{byr-cha-96}]. Byrne also did a further weakly nonlinear analysis
  and  numerical verification in [\ref{byr-99}].  Recently, Wu [\ref{wu-16}] proved that for
  $0<\tilde\sigma/\bar\sigma<\theta_*$ with some constant $\theta_*\in (0,1)$, problem (\ref{1.1}) has two
  radial stationary solutions, and there exist infinite many symmetry-breaking branches of stationary solutions
  bifurcating from radial stationary ones.

  If the boundary conditions $\sigma=\bar{\sigma}(1-\gamma\kappa)$ and
  $p=\bar{p}$ in (\ref{1.1}) are replaced
  by the following
\begin{equation}\label{1.2}
\sigma=\bar{\sigma},\quad  p=\gamma\kappa \quad \hbox{on $\partial\Omega(t) \times (0,+\infty)$},
\end{equation}
  the corresponding model is called {\em tumor spheroid model with surface tension effect}.
  Boundary condition  (\ref{1.2}) assumes that nutrient concentration is continuous across the tumor boundary
  and the pressure on the tumor surface is proportional to the mean curvature to maintain
  the cell-to-cell adhesiveness. Note that tumor spheroid model with surface tension
  effect belongs to category of Hele--Shaw type. In fact, by letting $\mu=0$, it becomes
  the classical Hele--Shaw problem.  The cell proliferation rate $\mu>0$ is indeed biologically
  meaningful in tumor spheroid model with surface tension effect. This novelty
  has attracted a lot of attention, and many illuminative results have been obtained.
  For $0<\tilde\sigma/\bar\sigma<1$, tumor spheroid model with surface tension effect
  has a unique radial stationary solution, which is globally asymptotically stable under
  radially symmetric perturbations (cf. [\ref{fri-rei-99}]). Bifurcation analysis shows that there
  exist infinite many symmetry-breaking branches of bifurcation stationary solutions (cf. [\ref{bor-fri-05},
  \ref{cui-esc-07}, \ref{fon-fri-03}, \ref{fri-rei-01}]). For asymptotic stability of radial stationary solution under
  non-radial perturbations, there exists a threshold value $\gamma_*$ of cell-to-cell adhesiveness
  such that the unique radial stationary solution is asymptotic stable for $\gamma>\gamma_*$, while it is unstable for
  $0<\gamma<\gamma_*$ (cf. [\ref{cui-09}, \ref{cui-esc-08}, \ref{fri-hu-06-0}, \ref{fri-hu-06}]).
  It is worthy to mention that there is not such a threshold value in stability result of Hele--Shaw problem.
  For extended studies of such type of tumor models, we refer readers to [\ref{cui-13}, \ref{fri-07}, \ref{fri-hu-07},
  \ref{fri-hu-08}, \ref{wu-cui-07}--\ref{wu-zhou-12},
  \ref{zhou-wu-15}] and references therein.

  In this paper, we study asymptotic stability of the two radial stationary solutions of (\ref{1.1})
  under non-radial perturbations.  Our motivation is two-fold. First, tumor spheroid model
  with surface tension effect neglects the discontinuity of nutrient flux across tumor boundary.
  Problem (\ref{1.1}) overcomes this disadvantage. As pointed out by Roose, Chapman and Maini,
  though the boundary condition with  Gibbs--Thomson relation seems speculative, it may be possible
  to check its veracity in experiments; moreover, problem (\ref{1.1}) raises a number of
  interesting mathematical and experimental points (see page 194 in [\ref{roose-07}]). Second,
  though the form of problem (\ref{1.1}) seems similar to tumor spheroid model
  with surface tension effect, the unknowns are coupled in different ways,
  which induces totally different nonlinearity and some new difficulties arise.
  In fact, based on classical theory of elliptic boundary value problem, problem (\ref{1.1}) and tumor spheroid
  model with surface tension effect can be both reduced into abstract differential equations of the form
  $\partial_t\rho+\Psi(\rho)=0$, where the unknown $\rho(t)$ is used to describe
  the free boundary $\partial\Omega(t)$. For tumor spheroid model with surface tension effect,
  $\Psi$ is a three-order quasi-linear differential operator in
  suitable Banach spaces. It can be regarded as a perturbation of the corresponding reduced operator of
  Hele--Shaw problem, so that its linearization is a sectorial operator and the corresponding abstract equation is of parabolic
  type. While for problem (\ref{1.1}), $\Psi$ is first-order and possesses a fully nonlinear structure.
  We have to make a delicate analysis of $\Psi$ and  determine whether
  its linearization is a sectorial operator or not.
   By using the technique of localization and Fourier multiplier analysis developed in [\ref{esc-sim-97-1}, \ref{esc-sim-97-2}],
  we overcome this difficulty and prove that the abstract equation of (\ref{1.1}) is still of parabolic type.
  Thus local well-posedness follows from analytic semigroup theory.
  Then with the help of some profound properties of modified Bessel functions, we calculate and analyze spectra of the
  linearized operator of $\Psi$. Finally,
  based on a {\em generalized principle of linearized stability} developed in [\ref{pru-sim-09}],
  we prove that the radial stationary solution with smaller radius is always unstable, and there exists a positive
  threshold value $\gamma_*$ of cell-to-cell adhesiveness $\gamma$, such that the radial stationary solution with
  larger radius is asymptotically stable for $\gamma>\gamma_*$, and unstable for $0<\gamma<\gamma_*$.

  It is interesting to compare tumor spheroid model with Gibbs--Thomson relation and tumor spheroid model with surface
  tension effect. For tumor spheroid model with Gibbs--Thomson relation, except for the two radial stationary solutions,
  the threshold value $\gamma_*$ of cell-to-cell adhesiveness, which is given by (\ref{4.25}) and (\ref{4.26}), is
  independent of proliferation rate $\mu$. This indicates that the proliferation rate $\mu$ has no effect on
  stability of radial dormant tumor. While for tumor spheroid model with surface tension effect, the corresponding
  threshold value $\gamma_*$ is a linear function of $\mu$ (cf. [\ref{cui-esc-08}, \ref{fri-hu-06}]), so that the
  proliferation rate $\mu$ indeed affects tumor's stability. This is a significant difference between these two
  modeling methods.

  To give a precise statement of our main result, we introduce some notations. Let $\Omega$ be a bounded open
  domain in $\Bbb R^3$ with a smooth boundary $\Gamma$. We denote by $BUC^{s}(\Omega)$ the space of all bounded
  and uniformly H\"older continuous functions on $\Omega$ of order $s>0$. Let $h^{s}(\Omega)$ denote the
  so-called little H\"older space, i.e., the closure of $BUC^\infty(\Omega)$ in $BUC^{s}(\Omega)$. Similarly,
  we denote by $h^{s}(\Gamma)$ the closure of $BUC^{\infty}(\Gamma)$ in $BUC^{s}(\Gamma)$.

  As mentioned before, if $0<\tilde\sigma/\bar\sigma<\theta_*$ then problem (\ref{1.1}) has two radial
  stationary solutions with radius denoted by $R_{s1}$ and $R_{s2}$, $(R_{s1}<R_{s2})$, respectively.
  Let  $R_s\in \{R_{s1}, R_{s2}\}$. Denote by $(\sigma_s,p_s,\Omega_s)$ the radial stationary solution
  with radius $R_s>0$, where $\Omega_s=\{x\in \Bbb R^3: |x|<R_s\}$, $\sigma_s$ and $p_s$ are given
  in (\ref{2.1}) and (\ref{2.2}). For a fixed $\alpha\in (0,1)$ and sufficiently small $\delta>0$, let
$$
  \mathcal O_\delta:=\{\rho\in h^{4+\alpha}(\Bbb S^2):  \|\rho\|_{h^{4+\alpha}(\Bbb S^2)}<\delta\}.
$$
  For given $\rho\in \mathcal O_\delta$, denote by $\Omega_\rho$ the domain enclosed by surface
  $\Gamma_\rho$, which is the image of the mapping $[\omega\to (R_s+\rho(\omega))\omega]$
  for $\omega \in \Bbb S^2$.  Thus solution $(\sigma,p,\Omega)$ can be rewritten as  $(\sigma, p, \rho)$
  with $\Omega=\Omega_\rho$, and $(\sigma_s,p_s,\Omega_s)$ can be rewritten as $(\sigma_s,p_s,
  \rho_s)$ with $\rho_s=0$.

  Let $T>0$, a triple
  $(\sigma, p, \rho)$ is called a classical solution of problem (\ref{1.1}) if
$$
  \sigma(\cdot,t), p(\cdot,t)\in h^{2+\alpha}(\Omega_{\rho(t)})\times h^{4+\alpha}(\Omega_{\rho(t)}),
   \; \; t\in (0,T],
$$
$$
  \rho(\cdot)\in C([0,T],\mathcal O_\delta) \cap C^1([0,T],h^{3+\alpha}(\Bbb S^2)),
$$
  and satisfies (\ref{1.1}) pointwise with $\Omega(t)=\Omega_{\rho(t)}$ and $\Omega_0=\Omega_{\rho_0}$.

  The first main result is on well-posedness theory of (\ref{1.1}), and is formulated below.
\medskip

   {\bf Theorem 1.1.} \ \ {\em  Given any initial data $\rho_0\in \mathcal O_\delta$
   for a sufficiently small $\delta>0$, there exists a unique classical solution
   $(\sigma,p,\rho)$ of problem $(\ref{1.1})$ on interval $[0,T]$ for some $T>0$.}
\medskip

  The proof of Theorem 1.1, given in Section 2 and Section 3,  is fulfilled by two steps. Firstly, by using
  Hanzawa transformation we reduce problem (\ref{1.1}) into an abstract differential equation
  $\partial_t\rho+\Psi(\rho)=0$, which contains an unknown function $\rho(t)$ describing the free boundary
  $\partial\Omega(t)$. Then we establish well-posedness by employing a localization technique, Fourier
  multiplier method and a delicate analysis of the linearization of $\Psi$.
\medskip

  Then we consider stability of each radial stationary solution $(\sigma_s,p_s,\Omega_s)$. An important character of
  problem (\ref{1.1}) is that it is invariant under coordinate translations.  For
  $(\sigma_s,p_s,\rho_s)$ with radius $R_s$ and a point $x_0\in \Bbb R^3$ in a small neighborhood of origin, we
  denote by $(\sigma^{[x_0]}_s,p_s^{[x_0]},\rho_s^{[x_0]})$ the translated radial stationary solution induced by
  $[x\to x+x_0]$.  Our second main result is stated as follows.
\medskip

 {\bf Theorem 1.2.} \ \  {\em  Let $0<\tilde\sigma/\bar\sigma<\theta_*$.
 There hold the following assertions:

  $(i)$  There exists a positive threshold value $\gamma_*$ of cell-to-cell adhesiveness such that for any $\gamma>\gamma_*$,
  the radial stationary solution $(\sigma_s, p_s, \rho_s)$ with larger radius $R_{s2}$ is asymptotically stable in the following
  sense: There exists constant $\epsilon\in(0,\delta)$ such that for any $\rho_0\in \mathcal O_\delta$ satisfying
  $\|\rho_0\|_{h^{4+\alpha}(\Bbb S^2)}<\epsilon$,  problem $(\ref{1.1})$ has a unique global classical
  solution $(\sigma,p,\rho)$, which converges exponentially fast to a translated radial stationary solution
  $(\sigma_s^{[x_0]},p_s^{[x_0]}, \rho_s^{[x_0]})$ in
  $h^{2+\alpha}(\Omega_{\rho(t)})\times h^{4+\alpha}
  (\Omega_{\rho(t)}) \times h^{4+\alpha}(\Bbb S^2)$
  for some $x_0\in \Bbb R^3 $, as $t\to+\infty$; while
  for $0<\gamma<\gamma_*$,
  it is unstable.

  $(ii)$ For all $\gamma>0$, the radial stationary solution $(\sigma_s, p_s, \rho_s)$
  with smaller radius $R_{s1}$ is unstable.}
\medskip

  Theorem 1.2 implies that cell-to-cell adhesiveness $\gamma$ plays a crucial role on
  tumor's stability and a smaller value of $\gamma$ may make tumor more aggressive.
  The proof of Theorem 1.2 is realized by first calculating all
  eigenvalues of the linearized operator at radial stationary solution
  and then employing the generalized
  principle of linearized stability developed in [\ref{pru-sim-09}].
\medskip

  The structure of the rest of this paper is arranged as follows. In the next section, by using Hanzawa
  transformation and classical theory of elliptic equations we reduce the free boundary problem (\ref{1.1})
  into an abstract differential equation in little H\"older spaces. Section 3 is devoted to establishing local
  well-posedness theory via localization argument and Fourier  multiplier method. In section 4 we study
  spectrum of linearized operator at radial stationary solution. In the last section
  we study asymptotic stability of each radial stationary solution and give the proof of Theorem 1.2.

\medskip
\hskip 1em

\section{Transformation and reduction}
\setcounter{equation}{0}
\hskip 1em

  In this section we transform free boundary problem (\ref{1.1})
  into an initial-boundary value problem on a fixed domain, and then
  reduce it into an evolutionary problem in Banach spaces.

  Firstly, we recall some results of problem ({1.1}) from [\ref{wu-16}].
  The radial stationary solution $(\sigma_s,p_s,\Omega_s)$
  can be given by
\begin{equation}\label{2.1}
 \Omega_s=\{x\in\Bbb R^3: r=|x|<R_s\},\qquad
 \sigma_s(r)=\bar\sigma(1-{\gamma\over R_s}){R_s\sinh r\over r\sinh R_s},
\end{equation}
\begin{equation}\label{2.2}
  p_s(r)=-\mu\bar\sigma(1-{\gamma\over R_s}){R_s\sinh r\over r\sinh R_s}
  +{1\over6}\mu\tilde\sigma r^2+\bar{p}+\mu\bar\sigma(1-{\gamma\over R_s})
  -{1\over6}\mu\tilde\sigma R_s^2,
\end{equation}
  where $R_s$ is a positive root of the following equation
\begin{equation}\label{2.3}
 f(R_s):= (1-{\gamma\over R_s}){R_s\coth R_s-1\over R_s^2}
 ={1\over 3}{\tilde\sigma\over\bar\sigma}.
\end{equation}
  The proof of Theorem 1.1 in [\ref{wu-16}] shows that there exists a constant
  $\theta_*\in(0,1)$ such that
\begin{equation}\label{2.4}
\begin{array}{c}
  \mbox{\em For } 0<\tilde\sigma/\bar\sigma<\theta_*,
  \mbox{ \em  equation
   $(\ref{2.3})$ has two positive roots } R_{s1} \mbox{ and } R_{s2}
   \\ [0.2 cm]
   \mbox{\em with }
   R_{s1}<R_{s2} \;\; and \;\; f'( R_{s1})>0,\,  f'(R_{s2})<0.
\end{array}
\end{equation}
  To economize our notation, in the sequel, we always set $R_s:=R_{si}$ for $i=1,2$.

   Denote $\omega(x):=x/|x|\in \mathbb S^2$ for $x\in \Bbb R^3\backslash \{0\}$.
   Taking $a_0>0$ sufficiently small, we see that the mapping
$$
  \Phi: \Bbb S^2\times(-a_0,a_0)\to\Bbb R^3,\quad
  \Phi(\omega,r)=(R_s+r) \omega,
$$
  is a $C^\infty$-diffeomorphism from $\Bbb S^2\times(-a_0,a_0)$
  onto its image $im(\Phi)$.

    Let $\delta\in (0, a_0/4)$ and $\alpha\in(0,1)$, set
$$
  \mathcal{O}_\delta:=\{\rho\in h^{4+\alpha}(\Bbb S^2):
  \,\|\rho\|_{h^{4+\alpha}(\Bbb S^2)}<\delta \}.
$$
  For each $\rho\in \mathcal{O}_\delta$ define
  $\Gamma_\rho:=im(\Phi(\cdot,\rho(\cdot))
  =\{(R_s+\rho(\omega))\omega:\omega\in\Bbb S^2\}$,
  and denote by $\Omega_\rho$ the domain enclosed by
  $\Gamma_\rho$.
  Clearly for the initial data $\Gamma_0=\partial\Omega_0$
  of class $h^{4+\alpha}$,
  there exists $\rho_0\in \mathcal{O}_\delta$ such that
  $\Gamma_{\rho_0}=\Gamma_0$,
  and accordingly, we have $\Omega_{\rho_0}=\Omega_0$.

  We take a cut-off function $\chi \in C^\infty(\Bbb R)$ such that
$$
  0\le\chi\le1,
\qquad
  \chi(\tau)=
  \left\{
  \begin{array}{l}
  1,\quad \mbox{for}\;|\tau| \le \delta,
  \\
  0,\quad \mbox{for}\;\;|\tau| \ge 3\delta,
  \end{array}
  \right.
  \qquad\mbox{and}\qquad
  |\chi'(\tau)|\le{2\over3\delta}.
$$
   Given $\rho\in\mathcal{O}_\delta$, we define
  the Hanzawa transformation
$$
  \Theta_\rho(x):=\left\{
  \begin{array}{ll}
  x+\chi(|x|-R_s)\rho(\omega(x))\omega(x) &\qquad\mbox{for}\;\; x\in im(\Phi),
  \\
  x& \qquad\mbox{for}\;\; x\notin im(\Phi).
 \end{array}
 \right.
$$
  Clearly one can verify that $\Theta_\rho(\Omega_s)=\Omega_\rho$,
  $\Theta_\rho(\Gamma_s)=\Gamma_\rho$, where $\Gamma_s=
  \partial\Omega_s=R_s\Bbb S^2$,
  and
$$
   \Theta_\rho\in \mbox{Diff}^{4+\alpha}(\Bbb R^3,\mathbb{R}^3)\cap
    \mbox{Diff}^{4+\alpha}(\Omega_s,\Omega_\rho).
$$
  Denote by $\Theta^\rho_*$ and $\Theta^*_\rho$
  respectively the push-forward and pull-back operators induced
  by $\Theta_\rho$, i.e.,
$$
  \Theta^\rho_* u=u\circ \Theta_\rho^{-1} \qquad \mbox{for}\;\;
  u\in BUC(\Omega_s),
$$
$$
 \;\Theta^*_\rho v=v \circ \Theta_\rho \qquad\; \mbox{  for}\;\;
  v\in BUC(\Omega_\rho).
$$
  By Lemma 2.1 in [\ref{esc-sim-97-2}], for $0\le k\le 4$ there holds
$$
  \Theta^*_\rho\in \mbox{Isom}(h^{k+\alpha}(\Omega_\rho),
  h^{k+\alpha}(\Omega_s)) \cap \mbox{Isom}(h^{k+\alpha}(\Gamma_\rho),
  h^{k+\alpha}(\Gamma_s))
$$
  with $\Theta^\rho_*=[\Theta^*_\rho]^{-1}$.
  Next, we define a function
  $\phi_\rho(x):=r(x)-R_s-\rho(\omega(x))$ for $x\in \Bbb R^3\backslash\{0\}$,
  where $r(x)=|x|$. Clearly,  $\Gamma_\rho=\phi_\rho^{-1}(0)$.
  We introduce the following transformed operators:
$$
  A(\rho)u:=\Theta_{\rho}^*\Delta(\Theta_*^{\rho}u),\qquad
  B(\rho)u=\Theta_\rho^*(\nabla\phi_\rho|_
  {\Gamma_\rho}\cdot\nabla(\Theta_*^\rho u)|_{\Gamma_\rho}),
$$
  for $u\in BUC^2(\Omega_s)$. Observe that $A(\rho)$ is the Laplace-Beltrami
  operator on $\Omega_s$ with respect to the metric induced by $\Theta_{\rho}$, and
  $B(\rho)$ is the induced outward normal derivative operator.
  Since $\Theta_\rho$ depends on $\rho$ analytically, by Lemma 2.2 in
   [\ref{esc-sim-97-1}], we have
\begin{equation}\label{2.5}
\left\{
\begin{array}{l}
  A(\cdot)\in C^\infty\big(\mathcal O_\delta,L(h^{4+\alpha}(\Omega_s),
 h^{2+\alpha}(\Omega_s))\big),
\\ [0.2 cm]
  B(\cdot) \in C^\infty\big(\mathcal O_\delta,L(h^{4+\alpha}(\Omega_s),
 h^{3+\alpha}(\Gamma_s))\big).
 \end{array}
 \right.
\end{equation}
  Let
$$
 \big( \mathcal B(\rho)u\big)(\omega)=\big(B(\rho)u\big)(x)
 \qquad\mbox{for}\;\;x\in \Gamma_s,\;\;\omega=x/|x|.
$$
  We have
\begin{equation}\label{2.6}
 \mathcal B(\cdot) \in C^\infty\big(\mathcal O_\delta,L(h^{4+\alpha}(\Omega_s),
 h^{3+\alpha}(\Bbb S^2))\big).
\end{equation}
  A direct computation shows that the mean curvature at point $x$ on $\Gamma_{\rho}$ is given by
\begin{equation}\label{2.7}
  \displaystyle\kappa(\rho)={1\over2}\left[{2r-\Delta_\omega\rho\over
  r\big(r^2+|\nabla_\omega \rho|^2\big)^{1/2}}+
  {2r|\nabla_\omega\rho|^2+\nabla_\omega|\nabla_\omega\rho|^2
  \cdot\nabla_\omega\rho\over  2r\big(r^2+|\nabla_\omega \rho|^2\big)^{3/2}}
  \right]_{r=R_s+\rho(\omega)},
\end{equation}
 where $\omega=x/|x|$ and $\Delta_\omega$ is  the Laplace-Beltrami operator
 on $\Bbb S^2$. We easily check that 
\begin{equation}\label{2.8}
  \kappa(\cdot)\in C^\infty(\mathcal O_\delta,h^{2+\alpha}(\Gamma_s)).
\end{equation}
  Finally, we set $u=\Theta^*_\rho\sigma$ and $v=\Theta^*_\rho p$.

  Using above notations,
  free boundary problem (\ref{1.1}) can be transformed into the following equivalent
  initial-boundary value problem
\begin{equation}\label{2.9}
\left\{
\begin{array}{ll}
    A(\rho)u=u \qquad \qquad  \qquad & \hbox{in  $\Omega_s
   \times [0,T]$},
   \\ [0.2 cm]
   A(\rho) v=-\mu(u-\tilde{\sigma})  \qquad\quad  &\hbox{in  $\Omega_s\times [0,T]$,}
\\ [0.2 cm]
  u=\bar{\sigma}(1-\gamma\kappa(\rho) ) \qquad\qquad & \hbox{on $\Gamma_s\times [0,T]$},
\\ [0.2 cm]
   v=\bar p \qquad \qquad  & \hbox{on $\Gamma_s\times [0,T]$},
\\ [0.2 cm]
     \partial_t \rho=-\mathcal B(\rho)v  \qquad \qquad & \hbox{on $\Bbb S^2\times [0,T]$},
\\ [0.2 cm]
   \rho(0)=\rho_0 \qquad \qquad  & \hbox{on $\Bbb S^2$}.
\end{array}
\right.
\end{equation}
   Here we used that the unit outward
  normal field ${\bf n}_{\rho(t)}$ and
  the normal velocity $V$ of $\Gamma_{\rho(t)}$
  are given by
$$
  {\bf n}_{\rho(t)}
  ={\nabla\phi_{\rho(t)}\over|\nabla\phi_{\rho(t)}|}
\qquad\mbox{and}\qquad
  V={\partial_t\rho\over |\nabla\phi_{\rho(t)}|},
$$
 respectively, for $\rho:=\rho(t)(\omega)\in C([0,T],\mathcal O_\delta)
 \cap C^1([0,T],h^{3+\alpha}(\Bbb S^2))$
  and some $T>0$.

  A triple $(u, v, \rho)$ is called a classical solution of problem (\ref{2.9}) if
$$
  (u, v)\in C([0,T],h^{2+\alpha}(\Omega_s)\times h^{4+\alpha}(\Omega_s)),
$$
$$
  \rho\in C([0,T],\mathcal O_\delta) \cap C^1([0,T],h^{3+\alpha}(\Bbb S^2)),
$$
  and it satisfies problem (\ref{2.9}) pointwise.
\medskip

  By above deduction we have

 \medskip

{\bf Lemma 2.1} \ \ {\em  A triple $(\sigma,p,\rho)$ is a classical
  solution of problem $(\ref{1.1})$ if and only if the triple $(u, v, \rho)$
  is a classical solution of problem $(\ref{2.9})$
 with $u=\Theta^*_\rho\sigma$ and $v=\Theta^*_\rho p$.
}

\medskip

  Next we further reduce problem (\ref{2.9}) into an abstract evolutionary problem containing free
 boundary $\rho$ only.  Given a function $\rho\in \mathcal O_\delta$,  a constant $\lambda\ge0$
 and $(g,h)\in  h^{k+\alpha}(\Omega_s)\times h^{k+2+\alpha}(\Gamma_s)$ for $k=0,1,2$,
 by well-known
 regularity theory of second-order elliptic equations, the elliptic boundary value problem
$$
\left\{
\begin{array}{ll}
  \lambda u- A(\rho) u=g\qquad & \mbox{in}\;\;\Omega_s,
  \\ [0.2 cm]
  u=h\qquad & \mbox{on}\;\; \Gamma_s,
\end{array}
\right.
$$
  has a unique solution which can be given by
\begin{equation}\label{2.10}
  u=\mathcal S_\lambda(\rho) g +\mathcal T_\lambda (\rho) h.
\end{equation}
  It holds that for $k=0,1,2$ (see Lemma 2.3 in [\ref{esc-sim-97-1}]),
\begin{equation}\label{2.11}
\left\{
\begin{array}{l}
 \mathcal S_\lambda\in C^\infty(\mathcal O_\delta, L(h^{k+\alpha}(\Omega_s),
  h^{k+2+\alpha}(\Omega_s)),
\\ [0.3 cm]
 \mathcal T_\lambda\in C^\infty(\mathcal O_\delta, L(h^{k+2+\alpha}(\Gamma_s),
  h^{k+2+\alpha}(\Omega_s)).
\end{array}
\right.
\end{equation}
  Thus for a given $\rho\in C([0,T],\mathcal O_\delta)\cap C^1([0,T],h^{3+\alpha}(\Bbb S^2))$,
  by solving the first four equations of problem (\ref{2.9}), we have
\begin{equation}\label{2.12}
\left\{
\begin{array}{l}
   u=\mathcal U(\rho):=\mathcal T_1(\rho) [\bar\sigma (1-\gamma\kappa(\rho))],
 \\ [0.2 cm]
   v=\mathcal V(\rho):=\mathcal S_0(\rho)\mathcal T_1(\rho) [ \mu\bar\sigma(1-\gamma\kappa(\rho))]
   -\mathcal S_0(\rho)\mu\tilde\sigma+ \bar p,
 \end{array}
\right.
\end{equation}
   where we  have used that  $\mathcal T_0(\rho)\bar p=\bar p$. Let
\begin{equation}\label{2.13}
  \Psi(\rho):=\mathcal B(\rho)\mathcal V(\rho)=-\mu\bar\sigma\gamma
  \mathcal B(\rho)\mathcal S_0(\rho)\mathcal T_1(\rho)\kappa(\rho)
  +\mathcal B(\rho)\mathcal S_0(\rho)\mathcal T_1(\rho)\mu\bar\sigma
  -\mathcal B(\rho)\mathcal S_0(\rho)\mu\tilde\sigma.
\end{equation}
  It follows from (\ref{2.5})--(\ref{2.8}) and (\ref{2.11}) that
\begin{equation}\label{2.14}
  \Psi  \in C^\infty(\mathcal O_\delta, h^{3+\alpha}(\Bbb S^2)).
\end{equation}
  By the above reduction, we see that problem (\ref{2.9}) is equivalent to
  the following
\begin{equation}\label{2.15}
\left\{
\begin{array}{ll}
  \partial_t \rho  + \Psi(\rho)=0
 & \qquad \mbox{on}\;\; \Bbb S^2\times [0,T],
  \\ [0.2 cm]
  \rho(0)=\rho_0  & \qquad \mbox{on}\;\; \Bbb S^2.
\end{array}
\right.
\end{equation}

 Similarly, a function $\rho\in C([0,T],\mathcal O_\delta) \cap
 C^1([0,T],h^{3+\alpha}(\Bbb S^2))$ is called a classical
 solution of problem (\ref{2.15}) if it satisfies each equation
 of problem (\ref{2.15}) on $[0,T]$ pointwise.

   In summary, we have

 \medskip

{\bf Lemma 2.2} \ \ {\em  The function $\rho$ is a classical
  solution of problem $(\ref{2.15})$ if and only if the triple
  $(u, v, \rho)$ is a classical solution of problem $(\ref{2.9})$
 with $(u,v)$ given by $(\ref{2.12})$.
}

\medskip

  In order to establish local well-posedness of the above problem (\ref{2.15}),
  we need to study the Fr\'echet derivative $\partial\Psi(\cdot)$ of nonlinear
  operator $\Psi$. In the sequel, we compute the Fr\'echet derivative
  $\partial \Psi(0)$ of $\Psi$ at $\rho=0$.

  Note that $\Omega_\rho\big|_{\rho=0}=\Omega_s$, and we have
$$
  \kappa(0)={1/ R_s},\qquad \mathcal U(0)=\sigma_s,
  \qquad
  \mathcal V(0)=p_s,\qquad
   \Psi(0)=p_s'(R_s)=0.
$$
  By (\ref{2.13}), $\Psi(\rho)=\mathcal B(\rho)\mathcal V(\rho)
  =\mu\mathcal B(\rho)\mathcal S(\rho)[\mathcal U(\rho)-\tilde\sigma]$.
  For any $\eta\in h^{4+\alpha}({\Bbb S^2})$,
\begin{equation}\label{2.16}
\begin{array}{rl}
  \partial \Psi(0)\eta\,\, &=\mu\mathcal B(0)\mathcal S_0(0)[\partial \mathcal U(0)\eta]
  +\mathcal B(0)\{\partial\mathcal S_0(0)[\eta,\mu(\sigma_s-\tilde\sigma)]\}
  +\partial \mathcal B(0)[\eta,p_s]
  \\ [0.2 cm]
  & =:\mbox{I}+\mbox{II}+\mbox{III},
\end{array}
\end{equation}
  where we use the notation
$$
  \displaystyle\partial \mathcal U(0)\eta:=\lim_{\varepsilon\to0}
  {\mathcal U(\varepsilon\eta)-\mathcal U(0)\over\varepsilon},
  \qquad
  \displaystyle\partial \mathcal B(0)[\eta,v]:=\lim_{\varepsilon\to0}
  {\mathcal B(\varepsilon\eta)v-\mathcal B(0)v\over
  \varepsilon}
$$
  for $\eta\in h^{4+\alpha}(\Bbb S^2)$ and
  $v\in h^{2+\alpha}(\Omega_s)$.

  Set
$$
  c_1:={\mu\gamma\bar\sigma\over2 R_s^2},\qquad
  c_2:={\mu\gamma\bar\sigma\over R_s^2}-\mu\sigma_s'(R_s),
  \qquad
  c_3:=p_s''(R_s).
$$

  We have  the following lemma:

 \medskip

 \medskip

{\bf Lemma 2.3} \ \ {\em  For $\eta\in h^{4+\alpha}(\Bbb S^2)$,
$$
  \partial\Psi(0)\eta=c_1\mathcal B(0)\mathcal S_0(0)\mathcal T_1(0)\Delta_\omega\eta
  +c_2\mathcal B(0)\mathcal S_0(0)\mathcal T_1(0)\eta+c_3\eta.
$$
}

    {\bf Proof}. \ \
  Since
   $\sigma_{\varepsilon\eta}:=\mathcal T_1(\varepsilon\eta)[\bar\sigma(1-\gamma\kappa(0))]$
   satisfies
\begin{equation}\label{2.17}
  \sigma_{\varepsilon\eta}-A(\varepsilon\eta)\sigma _{\varepsilon\eta} =0 \qquad \mbox{in}\;\;
  \Omega_s,
  \qquad\qquad \sigma_{\varepsilon\eta}=\bar\sigma(1-\gamma\kappa(0))\qquad\mbox{on}\;\;\Gamma_s,
\end{equation}
  and
  $\sigma_s=\mathcal T_1(0)[\bar\sigma(1-\gamma\kappa(0))]$ satisfies
\begin{equation}\label{2.18}
 \;\; \sigma_s-A(0)\sigma_s =0 \qquad \;\;\mbox{in}\;\;
  \Omega_s,
  \qquad\qquad \;\sigma_s=\bar\sigma(1-\gamma\kappa(0))
  \qquad\mbox{on}\;\;\Gamma_s,
\end{equation}
  we have
\begin{equation}\label{2.19}
  \big[\sigma_{\varepsilon\eta}-\sigma_s\big]
  -A(0)\big[\sigma_{\varepsilon\eta}-\sigma_s \big]
   -\big[ A(\varepsilon\eta)-A(0)\big]\sigma _{\varepsilon\eta}
  =0\qquad \mbox{in}\;\;
  \Omega_s.
\end{equation}
  Dividing both sides by $\varepsilon$ and letting $\varepsilon\to0$, we see
  that $\partial\mathcal T_1(0)[\eta,\bar\sigma(1-\gamma\kappa(0))]$ satisfies
\begin{equation}\label{2.20}
  u -A(0)u-\partial A(0)[\eta,\sigma_s]=0\qquad \mbox{in}\;\;
  \Omega_s,
   \qquad\qquad u=0 \qquad\mbox{on}\;\;\Gamma_s,
\end{equation}
  which implies that
\begin{equation}\label{2.21}
  \partial\mathcal T_1(0)[\eta,\bar\sigma(1-\gamma\kappa(0))]
  =\mathcal S_1(0)\{\partial A(0)[\eta,\sigma_s]\}.
\end{equation}
  Let $\sigma_{\varepsilon\eta}^s:=\Theta_{\varepsilon\eta}^*\sigma_s$,
  then
$$
u_0:=\lim_{\varepsilon\to0}{\sigma_{\varepsilon\eta}^s-\sigma_s\over\varepsilon}
=\chi(r-R_s)\sigma_s'(r)\eta.
$$
  Similarly as (\ref{2.20}) we find that $u_0$ satisfies
$$
  u_0 -A(0)u_0-\partial A(0)[\eta,\sigma_s]=0\qquad \mbox{in}\;\;
  \Omega_s,
   \qquad\qquad u_0=\sigma'(R_s)\eta \qquad\mbox{on}\;\;\Gamma_s,
$$
  and it follows that
$$
  u_0=\mathcal S_1(0)\{\partial A(0)[\eta,\sigma_s]\}+\mathcal T_1(0)(\sigma'(R_s)\eta).
$$
  Substituting it into (\ref{2.21}), we have
\begin{equation}\label{2.22}
\begin{array}{rl}
  \partial\mathcal T_1(0)[\eta,\bar\sigma(1-\gamma\kappa(0))]
  \,\,&=u_0-\mathcal T_1(0)(\sigma'(R_s)\eta)
  \\ [0.2 cm]
  &=\chi(r-R_s)\sigma_s'(r)\eta-\mathcal T_1(0)(\sigma'(R_s)\eta).
\end{array}
\end{equation}
  Let $p_{\varepsilon\eta}^s:=\Theta_{\varepsilon\eta}^*p_s$,
  then
$$
v_0:=\lim_{\varepsilon\to0}{p_{\varepsilon\eta}^s-p_s\over\varepsilon}
=\chi(r-R_s)p_s'(r)\eta.
$$
 Note that $p_s'(R_s)=0$, by a similar argument as above we can show that
\begin{equation}\label{2.23}
\begin{array}{rl}
  \partial\mathcal S_0(0)[\eta,\mu(\sigma_s-\tilde\sigma)]
  \,\,&=\mathcal S_0(0)\{\partial A(0)[\eta,p_s]\}
  \\ [0.2 cm]
  &=v_0-\mu\mathcal S_0(0)u_0.
\end{array}
\end{equation}
  By (\ref{2.7}), we easily have
\begin{equation}\label{2.24}
  \partial \kappa (0)\eta=\lim_{\varepsilon\to0}{\kappa(\varepsilon\eta)-\kappa(0)
  \over \varepsilon}= -{1\over R_s^2}(\eta+{1\over2}\Delta_\omega\eta).
\end{equation}
  Thus by (\ref{2.22}) and (\ref{2.24}) we compute
\begin{equation}\label{2.25}
\begin{array}{rl}
  \mbox{I} \,\,&= \mu\mathcal B(0)\mathcal S_0(0)[\partial\mathcal U(0)\eta]
  \\ [0.2 cm]
  &=  \mu\mathcal B(0)\mathcal S_0(0) \big\{\partial\mathcal T_1(0)
  [\eta,\bar\sigma(1-\gamma\kappa(0))]-\gamma\bar\sigma\mathcal T_1(0)\partial\kappa(0)\eta
  \big\}
  \\ [0.2 cm]
  &= \mu\mathcal B(0)\mathcal S_0(0) \big\{u_0-\mathcal T_1(0)(\sigma_s'(R_s)\eta)
  -\gamma\bar\sigma\mathcal T_1(0)\partial\kappa(0)\eta
  \big\}
  \\ [0.2 cm]
  &\displaystyle=\mu\mathcal B(0)\mathcal S_0(0)u_0 +
  {\mu\gamma\bar\sigma\over2 R_s^2}
   \mathcal B(0)\mathcal S_0(0)\mathcal T_1(0)\Delta_\omega\eta+
 \mu( {\gamma\bar\sigma\over R_s^2}-\sigma_s'(R_s))
 \mathcal B(0)\mathcal S_0(0)\mathcal T_1(0)\eta.
\end{array}
\end{equation}
  By (\ref{2.23}) we have
\begin{equation}\label{2.26}
  \begin{array}{rl}
 \mbox{II} \,\,&= \mathcal B(0)\big\{
  \partial\mathcal S_0(0)[\eta, \mu(\sigma_s-\tilde\sigma)]
  \big\}
  \\ [0.2 cm]
  &=\mathcal B(0)v_0-\mu\mathcal B(0)\mathcal S_0(0)u_0
  \\ [0.2 cm]
  &= p_s''(R_s)\eta-\mu\mathcal B(0)\mathcal S_0(0)u_0.
\end{array}
\end{equation}
  Finally, a direct computation shows that
\begin{equation}\label{2.27}
  \mbox{III}=\lim_{\varepsilon\to0}{\mathcal B(\varepsilon\eta)p_s
  -\mathcal B(0)p_s\over \varepsilon}
  =p_s'(R_s)\eta=0.
\end{equation}
  Hence, by (\ref{2.16}) and adding (\ref{2.25})--(\ref{2.27}),
  we immediately complete the proof.
  \qquad$\Box$

\medskip
\hskip 1em
\section{Local well-posedness}
\setcounter{equation}{0}
\hskip 1em

  In this section, we prove local well-posedness of problem (\ref{2.15}).
  By using a technique of localization and Fourier multiplier method
  developed in [\ref{esc-sim-97-1}, \ref{esc-sim-97-2}], we shall prove
  $-\partial\Psi(0)$ is an infinitesimal generator of a strongly continuous
  analytic semigroup on $h^{3+\alpha}(\Bbb S^2)$, and
  problem (\ref{2.15}) is of parabolic type for $\delta$ is sufficiently small.
  Thus by using analytic semigroup theory and applications to parabolic
  differential problems (cf. [\ref{amann}, \ref{lunardi}]), we get local
  well-posedness.

  For giving expression of linear operators introduced in above section
  by means of local coordinates, we need some notations similarly
  as in [\ref{esc-gui}, \ref{esc-sim-97-1}, \ref{esc-sim-97-2}].   Given $a\in (0,a_0)$.
  We set $\mathcal R_a:=\{x\in \Bbb R^3: R_s-a<|x|\le R_s\}$ which
  is a neighborhood of $\Gamma_s$.
  There exist an integer $m:=m_a\in \Bbb N$ and an atlas $\{(U_l,\psi_l); 1\le l\le m\}$\
  of $\mathcal R_a$ such that
$$
   \psi_l\in C^{\infty}(U_l,(-a,a)^2\times [0,a)), \qquad
   \psi_l:U_l\cap \Gamma_s\to (-a,a)^2\times \{0\},\qquad
   1\le l\le m.
$$
   Given $l\in\{1,\cdots,m\}$, denote by $\psi_*^l$ and $\psi^*_l$
   the push-forward and pull-back operators
   induced by $\psi_l$, respectively, i.e.,
$$
  \psi^l_* u=u\circ \psi_l^{-1} \quad\;\; \mbox{for}\;\;
  u\in BUC(U_l),\qquad\qquad\quad\;\;
$$
$$
 \psi^*_l v=v \circ \psi_l \qquad\; \mbox{  for}\;\;
  v\in BUC((-a,a)^2\times [0,a)).
$$
  Define the local representation $\tilde {\mathcal A}_l$,
  $\tilde {\mathcal B}_l$
  and $\tilde { \mathcal P}_l$ of $A(0)$, $\mathcal B(0)$
  and $\Delta_\omega$,
  respectively, on $((-a,a)^2\times [0,a),\psi_l)$ as
$$
  \tilde{\mathcal A}_l:=\psi^l_*\circ A(0)\circ \psi^*_l ,\qquad
  \tilde{\mathcal B}_l:=\psi^l_*\circ \mathcal B(0)\circ \psi^*_l ,\qquad
  \tilde{\mathcal P}_l:=\psi^l_*\circ \Delta_\omega\circ \psi^*_l.
$$
  Since $A(0)$ is the Laplace operator on $\Omega_s$,
  $\mathcal B(0)$ is the outward normal derivative operator, and
  $\Delta_\omega$ is the Laplace-Beltrami operator on $\Bbb S^2$,
  there exist (cf. Lemma 3.2 in [\ref{esc-sim-97-1}] and Lemma 3.1 in
  [\ref{esc-sim-97-2}])
$$
  a_{jk}^l(\cdot), a_j^l(\cdot) \in C^\infty((-a,a)^2\times [0,a),\Bbb R),\quad
   b_j^l(\cdot)\in C^\infty((-a,a)^2,\Bbb R),
  \qquad 1\le j,k\le 3,
$$
$$
  p_{jk}^l(\cdot), p_j^l(\cdot)\in C^\infty((-a,a)^2,\Bbb R),
  \qquad 1\le j,k\le 2,
$$
 such that
\begin{equation}\label{3.1}
\begin{array}{l}
  [a_{jk}^l] \mbox{ is symmetric and uniformly positive definite on }
  (-a,a)^2\times [0,a),
\\ [0.3 cm]
 [p_{jk}^l] \mbox{ is symmetric and uniformly positive definite on }
 (-a,a)^2,
\\ [0.3 cm]
  b_3^l \mbox{ is uniformly positive on }
  (-a,a)^2,
\end{array}
\end{equation}
  and such that
\begin{equation}\label{3.2}
\begin{array}{c}
 \displaystyle\tilde{\mathcal A}_l = \sum_{j,k=1}^3 a_{jk}^l(\cdot)\partial_j\partial_k
  +\sum_{j=1}^3 a_j^l(\cdot) \partial_j,
  \qquad
  \displaystyle \tilde{\mathcal B}_l = -\sum_{j=1}^3 b_j^l(\cdot)
   \Upsilon\partial_j,
  \\ [0.3 cm]
    \displaystyle  \tilde{\mathcal P}_l = \sum_{j,k=1}^2 p_{jk}^l(\cdot)\partial_j\partial_k
   +\sum_{j=1}^2 p_j^l(\cdot) \partial_j,
\end{array}
\end{equation}
  where $\Upsilon$ is the trace operator on $\Bbb R^2\times \{0\}$.

  Next we fix the localization at $x_l:= \psi_l^{-1}(0)\in U_l\cap \Gamma_s$
  and define the following linear differential
  operators with constant coefficients:
$$
 \displaystyle \tilde{ \mathcal A}_{l,0}:=\sum_{j,k=1}^3 a_{jk}^l(0)
 \partial_j\partial_k,
  \qquad
  \displaystyle \tilde{\mathcal B}_{l,0} := -\sum_{j=1}^3 b_j^l(0)
   \Upsilon\partial_j,
  \qquad
    \displaystyle  \tilde{\mathcal P}_{l,0} := -1+\sum_{j,k=1}^2 p_{jk}^l(0)\partial_j\partial_k.
$$
  For any given $(g,h)\in h^{k+\alpha}(\Bbb H^3)\times h^{k+2+\alpha}(\partial \Bbb H^3)$,
  $k=0,1,2$, we consider the following elliptic boundary value problem in the half space
   $\Bbb H^3:=\Bbb R^2\times (0,+\infty)$:
\begin{equation}\label{3.3}
\left\{
\begin{array}{ll}
  u-\tilde{\mathcal A}_{l,0} u=g\qquad & \mbox{in}\;\; \Bbb H^3,
  \\ [0.2 cm]
  \Upsilon u=h\qquad & \mbox{on}\;\; \partial\Bbb H^3:=\Bbb R^2\times\{0\}.
\end{array}
\right.
\end{equation}
  By Appendix B in [\ref{esc-sim-95}],  there exists a unique solution
  $u\in h^{k+2+\alpha}(\Bbb H^3)$ of problem (\ref{3.3}) which is given by
  $ u=\tilde{\mathcal S}_{l,0} g+\tilde{\mathcal T}_{l,0} h$ and
$$
  \tilde{\mathcal S}_{l,0}\in L(h^{k+\alpha}(\Bbb H^3),h^{k+2+\alpha}(\Bbb H^3)),
$$
$$
  \tilde{\mathcal T}_{l,0}\in L(h^{k+2+\alpha}(\partial\Bbb H^3),h^{k+2+\alpha}(\Bbb H^3)),
$$
   for $k=0,1,2$. Moreover,  $\tilde{\mathcal S}_{l,0}$ and $\tilde{\mathcal T}_{l,0}$
   can be expressed as Fourier multiplier operators. For this purpose, we let
$$
  \mathbf a^l:=(a_{13}^l(0), a_{23}^l(0)),\qquad
  a^l_0(\xi):=\sum_{j,k=1}^2a_{jk}^l(0)\xi_1\xi_2,
$$
$$
  \lambda^l(\xi):={i\langle\mathbf a^l, \xi\rangle\over a_{33}^l(0)}
  +{1\over a_{33}^l(0)}\sqrt{a_{33}^l(0)(1+a^l_0(\xi))-
  \langle\mathbf a^l, \xi\rangle^2},
$$
  and
$$
  q^l(\xi,z):=1+a^l_0(\xi)+2i\langle\mathbf a^l, \xi\rangle z-a_{33}^l(0)z^2,
  \qquad \xi=(\xi_1,\xi_2)\in\Bbb R^2,
  \;\;z\in\Bbb C,
$$
   where $\langle\cdot,\cdot\rangle$ denotes the inner product in $\Bbb R^2$.
   Due to (\ref{3.1}), we easily verify that $\lambda^l(\xi)$ is well-defined and for
   any given $\xi\in\Bbb R^2 $, $\lambda^l(\xi)$ is the unique root with positive
   real part of $q^l(\xi,z)$.
   Denote $\mathcal F$ and $\mathcal F^{-1}$ by the Fourier transform and
   inverse Fourier transform in $\Bbb R^2$, respectively. We have
   (cf. (4.2)--(4.4) in [\ref{esc-sim-97-2}])
\begin{equation}\label{3.4}
   \tilde{\mathcal S}_{l,0} g(x,y)=\Big[\mathcal F^{-1} (1-e^{-\lambda^l(\cdot)y})
   {1\over 1+a_0^l(\cdot)} \mathcal F g\Big](x),
\end{equation}
\begin{equation}\label{3.5}
   \tilde{\mathcal T}_{l,0} h(x,y)=
  \Big[\mathcal F^{-1} e^{-\lambda^l(\cdot)y} \mathcal F h\Big](x),
\end{equation}
  for $g\in h^{k+2+\alpha}(\Bbb H^3)$,
 $h\in h^{k+\alpha}(\Bbb R^2)$
 and $(x,y)\in \Bbb H^3$, $k=0,1,2$.
  Let
$$
  p^l_0(\xi):=\sum_{j,k=1}^2 p_{jk}^l(0)\xi_j\xi_k
  \qquad \mbox{for}\;\; \xi\in\Bbb R^2.
$$
  It is easy to see that
\begin{equation}\label{3.6}
  \tilde{\mathcal P}_{l,0}h (x)=-\Big[\mathcal F^{-1}(1+ p_0^l(\cdot))
  \mathcal F h\Big](x)
  \qquad \mbox{for}\;\; h^{k+2+\alpha}(\Bbb R^2),
\end{equation}
  and
$$
  \tilde{\mathcal P}_{l,0}\in L(h^{k+2+\alpha}(\Bbb R^2),h^{k+\alpha}(\Bbb R^2))
  \qquad\mbox{for}\;\; k=0,1,2.
$$
  By (\ref{3.4})--(\ref{3.6}), we have
\begin{equation}\label{3.7}
  \tilde{\mathcal S}_{l,0}\tilde{\mathcal T}_{l,0}\tilde{\mathcal P}_{l,0}h(x,y)=
  \Big[\mathcal F^{-1} \Big\{
  (e^{-\lambda^l(\cdot)y}-1) e^{-\lambda^l(\cdot)y}
   {1+ p_0^l(\cdot)\over 1+a_0^l(\cdot)}
    \Big\} \mathcal F h\Big](x),
\end{equation}
  for $h\in h^{4+\alpha}(\Bbb R^2)$ and $(x,y)\in \Bbb H^3$.
  Define
$$
  g^l(\xi):=b_3^l(0)\lambda^l(\xi){1+ p_0^l(\xi)\over 1+a_0^l(\xi)}\qquad
  \mbox{for}\;\;\xi\in\Bbb R^2.
$$
  By (\ref{3.7}) we immediately have
\begin{equation}\label{3.8}
  \mathcal G_l:= \tilde{\mathcal B}_{l,0}\tilde{\mathcal S}_{l,0}\tilde{\mathcal T}_{l,0}
  \tilde{\mathcal P}_{l,0}
  =\mathcal F^{-1}  g^l(\cdot) \mathcal F .
\end{equation}

    Given two Banach spaces $E_0$ and $E_1$ such that $E_1$
    is densely and continuously embedded into $E_0$,  we denote
    $\mathcal H(E_1,E_0)$ by the subspace of all linear operators
    $A\in L(E_1,E_0)$ such that $-A$ generates a strongly continuous
    analytic semigroup on $E_0$. By Remark I.1.2.1 (a) in [\ref{amann}],
    we have that a linear operator $A\in L(E_1,E_0)$ belongs to
    $\mathcal H(E_1,E_0)$ if and only if there exist positive
    constants $C$ and $\lambda_*$ such that
 \begin{equation}\label{3.9}
 \begin{array}{l}
    \lambda_*+A\in {\rm Isom}(E_1,E_0),
    \\ [0.3 cm]
    |\lambda|\|x\|_{E_0}+\|x\|_{E_1}
    \le C\|(\lambda+A)x\|_{E_0},
\end{array}
\end{equation}
    for all $x\in E_1$ and $ {\rm Re}\lambda\ge \lambda_*$.

    Based on the property of the symbol $g_l(\xi)$ of $\mathcal G_l$ and
    Fourier multiplier theory, we have

\medskip

\medskip

{\bf Lemma 3.1}\ \ $\mathcal G_l\in \mathcal H(h^{4+\alpha}(\Bbb R^2),
   h^{3+\alpha}(\Bbb R^2))$.

\medskip

 {\bf Proof}. \ \  Let
$$
  {\tilde\lambda}(\xi,\tau):={i\langle\mathbf a^l, \xi\rangle\over a_{33}^l(0)}
  +{1\over a_{33}^l(0)}\sqrt{a_{33}^l(0)(\tau^2+a^l_0(\xi))-
  \langle\mathbf a^l, \xi\rangle^2},
$$
  and
$$
  {\tilde g}(\xi,\tau):= b_3^l(0)\tilde \lambda(\xi,\tau){\tau^2+ p_0^l(\xi)\over \tau^2
  +a_0^l(\xi)}\qquad \mbox{for}\;\;\xi\in\Bbb R^2,\;\;
  \tau\in (0,+\infty).
$$
  It is easy to verify that
$$
  \lambda^l(\xi)={\tilde\lambda}(\xi,1)\qquad
  \mbox{and}\qquad
  g^l(\xi)=\tilde g(\xi,1).
$$
  Moreover,
  $\tilde g(\xi,\tau)$ is positively homogeneous
  of degree 1 and all derivatives of $\tilde g(\xi,\tau)$ are bounded on
  $|\xi|^2+\tau^2=1$.

   By (\ref{3.1}), we see that $b_3^l(0)>0$ and $[a_{jk}^l(0)]$,
   $[p_{jk}^l(0)]$ are symmetric and positive definite. It implies that
   there exist positive constants $a_*$, $b_*$, $c_*$ such that
   for all $\xi\in\Bbb R^2$ and $\tau\in(0,+\infty)$,
$$
  {\rm Re} \tilde \lambda(\xi,\tau)\ge a_*\sqrt{|\xi|^2+\tau^2},
\quad
  \tau^2+a_0^l(\xi)\le b_*(|\xi|^2+\tau^2),
\quad
  \tau^2+p_0^l(\xi)\ge c_*(|\xi|^2+\tau^2).
$$
  Thus we have
$$
\begin{array}{rl}
  {\rm Re} \tilde g(\xi,\tau)\,\,&=\displaystyle b_3^l(0)
  {\rm Re}\tilde \lambda(\xi,\tau)
  {\tau^2+ p_0^l(\xi)\over \tau^2+a_0^l(\xi)}
  \\ [0.3 cm]
 & \ge d_* \sqrt{\xi^2+\tau^2},
\end{array}
$$
  where $d_*:=a_*b_*^{-1}c_*b_3^l(0)>0$. With these properties of
  $\tilde g(\xi,\tau)$ and by
  using Mikhlin-H\"ormander multiplier theorem, we have
  $\mathcal G_l=\mathcal F^{-1}  \tilde g(\cdot,1) \mathcal F
  \in \mathcal H(h^{4+\alpha}(\Bbb R^2),
   h^{3+\alpha}(\Bbb R^2))$,
  cf. the proof of Theorem 4.2 in [\ref{esc-sim-97-2}] or
  Theorem A.2 in [\ref{esc-sim-95}]. The proof is complete.
 \qquad$\Box$

\medskip

 Introduce an operator $\vartheta_*: BUC (\Bbb S^2)\to BUC(\Gamma_s)$
 given by
$$
  \vartheta_*\eta(R_s\omega)=\eta(\omega)\qquad
  \mbox{for}\;\;\eta\in BUC(\Bbb S^2),\;\;\omega\in\Bbb S^2.
$$
  Let $\beta\in(0,\alpha)$ be fixed and
$$
  G:=\mathcal B(0)\mathcal S_0(0)\mathcal T_1(0)\Delta_\omega.
$$
  For any $\epsilon>0$, by
 the proof of Lemma 5.1 in [\ref{esc-sim-97-2}] (with a trivial modification),
 there exist $a\in (0,a_0]$ and a partition of unity $\{(U_l,\psi_l);1\le l\le m_a\}$
 for $\mathcal R_a$, and a positive constant $C:=C(\epsilon,a)$ such that
\begin{equation}\label{3.10}
  \|\psi_*^l(\psi_l \vartheta_* G \eta)-\mathcal G_l \psi_*^l(\psi_l\vartheta_*\eta)
  \|_{h^{3+\alpha}(\Bbb R^2)}
  \le \epsilon \|\psi_*^l(\psi_l\vartheta_*\eta)\|_{h^{4+\alpha}(\Bbb R^2)}+
  C\|\eta\|_{h^{4+\beta}(\Bbb S^2)},
 \end{equation}
  for all $\eta\in h^{4+\alpha}(\Bbb S^2)$ and $l=1,2,\cdots,m_a $.

  With the above preparation, we have the following

\medskip

\medskip

{\bf Lemma 3.2} \ \ $\partial \Psi(0)\in
  \mathcal H(h^{4+\alpha}(\Bbb S^2),
   h^{3+\alpha}(\Bbb S^2))$.

\medskip

 {\bf Proof}. \ \ $(i)$  Based on Lemma 3.1 and the sharp
   estimate (\ref{3.10}), by applying interpolation inequality
   and (\ref{3.9}) we can show that
\begin{equation}\label{3.11}
  G=\mathcal B(0)\mathcal S_0(0)\mathcal T_1(0)\Delta_\omega
  \in \mathcal H(h^{4+\alpha}(\Bbb S^2),  h^{3+\alpha}(\Bbb S^2)).
\end{equation}
   We omit the details and refer to the proof of Theorem 5.2 in
   [\ref{esc-sim-97-2}] with a minor modification, see also a similar
   proof of Theorem 4.1 in [\ref{esc-sim-97-1}].

  $(ii)$ Recall that $c_1={1\over 2 R_s^2}\mu\gamma\bar\sigma>0$, by
  (\ref{3.11}) we have $
  c_1G\in \mathcal H(h^{4+\alpha}(\Bbb S^2),  h^{3+\alpha}(\Bbb S^2))$.
  Let
$$
  K\eta:= c_2 \mathcal B(0)\mathcal S_0(0)\mathcal T_1(0)\eta+c_3  \eta
  \qquad\mbox{for}\;\;\eta\in h^{4+\alpha}(\Bbb S^2).
$$
  By (\ref{2.6}) and  (\ref{2.11}), we see that
$$
 K \in L(h^{4+\alpha}(\Bbb S^2),  h^{4+\alpha}(\Bbb S^2)).
$$
  Since $h^{4+\alpha}(\Bbb S^2)$ is compactly embedded into
  $h^{3+\alpha}(\Bbb S^2)$, by a standard perturbation result
  (cf. Theorem I.1.5.1 in [\ref{amann}]),
  we have $\partial \Psi(0)=c_1G+K\in  \mathcal H(h^{4+\alpha}(\Bbb S^2),
  h^{3+\alpha}(\Bbb S^2))$. The proof is complete.
  \qquad $\Box$

\medskip

  Since $\mathcal H(h^{4+\alpha}(\Bbb S^2),
  h^{3+\alpha}(\Bbb S^2))$ is an open subset of
  $L(h^{4+\alpha}(\Bbb S^2), h^{3+\alpha}(\Bbb S^2))$
  (cf. Theorem I.1.3.1 in [\ref{amann}]), by Lemma 3.2 we
  immediately get

\medskip

\medskip

{\bf Corollary 3.3}\ \ {\em Let $\delta>0$ be sufficiently small.
  For any given $\rho \in \mathcal O_\delta$, there holds
$$
  \partial\Psi(\rho)\in \mathcal H(h^{4+\alpha}(\Bbb S^2),
  h^{3+\alpha}(\Bbb S^2)).
$$
}

  Corollary 3.3 implies that problem (\ref{2.15}) in $\mathcal O_\delta$
  is of parabolic type in the sense of Amann [\ref{amann}] and
   Lunardi [\ref{lunardi}]. By applying analytic semigroup theory and
   applications to parabolic differential problems (cf. Theorem 8.1.1
   and Theorem 8.3.4 in [\ref{lunardi}]), we get the following
   result:

\medskip

\medskip

{\bf Theorem 3.4}\ \  {\em Given $\rho_0\in \mathcal O_\delta$. There exists
  $T>0$ such that problem $(\ref{2.15})$ has a unique classical solution
  $\rho\in C([0,T],\mathcal O_\delta)\cap C^1([0,T], h^{3+\alpha}(\Bbb S^2))$.
 }

\medskip

  By Lemma 2.1, Lemma 2.2 and Theorem 3.4, we see that Theorem 1.1 follows
  and local well-posedness of free boundary problem (\ref{1.1}) has been obtained.

\medskip

\hskip 1em

\section{Spectrum analysis}
\setcounter{equation}{0}
\hskip 1em

  In this section we study the spectrum of $\partial\Psi(0)$.
  Since $h^{4+\alpha}(\Bbb S^2)$ is compactly embedded
  into $h^{3+\alpha}(\Bbb S^2)$, by Lemma 3.2, we see that
  the spectrum of $\partial\Psi(0)$ consists of all eigenvalues,
  which will be obtained by employing spherical harmonics
  and modified Bessel functions.

  Recall that $\Delta_\omega$ is denoted by the
  Laplace-Beltrami operator on the unit sphere $\Bbb S^2$, and
$$
  \Delta={\partial^2\over\partial r^2}+{2\over r}{\partial
  \over\partial r}+{1\over r^2}\Delta_\omega.
$$
  By Lemma 2.3, we have for given $\eta\in h^{4+\alpha}(\Bbb S^2)$,
$$
  \partial\Psi(0)\eta=c_1\mathcal B(0)\mathcal S_0(0)\mathcal T_1(0)\Delta_\omega\eta
  +c_2\mathcal B(0)\mathcal S_0(0)\mathcal T_1(0)\eta+c_3\eta,
$$
  where
$$
  c_1={\mu\gamma\bar\sigma\over2 R_s^2},\qquad
  c_2={\mu\gamma\bar\sigma\over R_s^2}-\mu\sigma_s'(R_s),
  \qquad
  c_3=p_s''(R_s).
$$
  By using (\ref{2.1})--(\ref{2.3}), a direct calculation shows that
$$
  \sigma_s'(R_s)={1\over 3}\tilde\sigma R_s\qquad
  \mbox{and}\qquad
  p_s''(R_s)=-\mu\big[\bar\sigma(1-{\gamma\over R_s})
  -\tilde\sigma\big].
$$
  Hence we can rewrite
\begin{equation}\label{4.1}
  c_1={\mu\gamma\bar\sigma\over2 R_s^2},\qquad
  c_2={\mu\gamma\bar\sigma\over R_s^2}-{1\over 3}\mu\tilde\sigma R_s,
  \qquad
  c_3=-\mu\big[\bar\sigma(1-{\gamma\over R_s})
  -\tilde\sigma\big].
\end{equation}
  Consider the following problem
\begin{equation}\label{4.2}
\left\{
\begin{array}{ll}
   u- \Delta u=0
 \qquad & \mbox{in}\;\;  \Omega_s,
\\ [0.3 cm]
  u =c_1\Delta_\omega \eta +c_2\eta
  \qquad\qquad & \mbox{on}\;\;  \Gamma_s,
\\ [0.3 cm]
   -\Delta v= u
    &  \mbox{in}\;\;  \Omega_s,
\\ [0.3 cm]
  v=0   & \mbox{on}\;\;  \Gamma_s,
\end{array}
\right.
\end{equation}
  where $u=u(r,\omega)$ and $v=v(r,\omega)$ are unknown functions.
  By (\ref{2.10}), we see that the solution of problem (\ref{4.2}) is given
  by
\begin{equation}\label{4.3}
  u=\mathcal T_1(0)[c_1\Delta_\omega \eta +c_2\eta],
  \qquad
  v=\mathcal S_0(0)\mathcal T_1(0)[c_1\Delta_\omega \eta +c_2\eta].
\end{equation}
  Note that $\displaystyle \mathcal B(0)v={\partial v\over \partial r}\Big|_{r=R_s}$
  for $v\in BUC^2(\Omega_s)$,
  we have
\begin{equation}\label{4.4}
\begin{array}{rl}
  \partial\Psi(0)\eta \,\,&= \mathcal B(0)v+c_3\eta
\\ [0.3 cm]
  &\displaystyle={\partial v\over \partial r}\Big|_{r=R_s}
  -\mu\big[\bar\sigma(1-{\gamma\over R_s})
  -\tilde\sigma\big]\eta,
\end{array}
\end{equation}
  where $v=v(r,\omega)$ is given by (\ref{4.3}).

  Next we solve problem (\ref{4.2}) for given
\begin{equation}\label{4.5}
   \eta(\omega)=\sum_{k=0}^\infty\sum_{l=-k}^{k}
   c_{kl}Y_{k,l}(\omega)\in C^\infty(\Bbb S^2),
\end{equation}
  where $Y_{k,l}(\omega)$ $(k\ge0,-k\le l\le k)$ denotes
  spherical harmonic of order $(k,l)$, and
  $c_{kl}$ is rapidly decreasing in $k$.
  By regularity theory of elliptic differential equations, we see that
  there is a unique smooth solution $(u,v)$ of problem (\ref{4.2}).
  Let
\begin{equation}\label{4.6}
   \displaystyle u(r,\omega)=\sum_{k=0}^\infty\sum_{l=-k}^{k}
     a_{kl}(r)Y_{k,l}(\omega)\qquad \mbox{and}
\qquad  v(r,\omega)=\sum_{k=0}^\infty\sum_{l=-k}^k
     b_{kl}(r)Y_{k,l}(\omega),
\end{equation}
  where $a_{kl}(r)$ and $b_{kl}(r)$ are unknown functions
  to be determined later.
   Recall that the well-known formula
$$
  \Delta_\omega Y_{k,l}(\omega)=-(k^2+k) Y_{k,l}(\omega)\qquad\mbox{for}\;\;
  k\ge0.
$$
  By substituting (\ref{4.5}) and (\ref{4.6}) into (\ref{4.2}), and
  comparing coefficients of each $Y_{k,l}(\omega)$, we have
\begin{eqnarray}\label{4.7}
  &&a''_{kl}(r)+{2\over r} a'_{kl}(r)-{k^2+k\over
  r^2}a_{kl}(r)=a_{kl}(r),
\\ [0.3cm]\label{4.8}
  &&a_{kl}(R_s)=\big[c_2-c_1(k^2+k)\big]c_{kl},
\\ [0.3cm]\label{4.9}
  &&b''_{kl}(r)+{2\over r} b'_{kl}(r)-{k^2+k\over
  r^2}b_{kl}(r)=- a_{kl}(r),
\\ [0.3cm]\label{4.10}
  &&b_{kl}(R_s)=0.
\end{eqnarray}
  Recall modified Bessel functions (cf. [\ref{wat-44}]):
\begin{equation}\label{4.11}
  I_m(r)=\sum_{k=0}^\infty {(r/2)^{m+2k}\over k!\Gamma(m+k+1)}
  \qquad \mbox{for}\;\; m\ge 0,
\end{equation}
  which satisfies
\begin{equation}\label{4.12}
\left\{
\begin{array}{l}
  \displaystyle I''_m(r)+{1\over r}I'_m(r)-(1+{m^2\over r^2})\,I_m(r)=0
  \qquad \mbox{for} \;\; r>0,
  \\ [0.3cm]
  I_m(r) \;\;\mbox{bounded at}\;\; r\sim 0,
\end{array}
\right.
\end{equation}
  and there hold the following properties:
\begin{equation}\label{4.13}
  I'_m(r)+{m\over r}I_m(r)=I_{m-1}(r) \qquad \mbox{for}\;\;m\ge1,
\end{equation}
\begin{equation}\label{4.14}
  I'_m(r)-{m\over r}I_m(r)=I_{m+1}(r) \qquad \mbox{for}\;\;m\ge0.
\end{equation}
   By using (\ref{4.12}) we get the solution of problem (\ref{4.7})--(\ref{4.10})
   is given by
\begin{equation}\label{4.15}
   a_{kl}(r)=c_{kl}\big[c_2-c_1(k^2+k)\big]
   {R_s^{1\over2} I_{k+{1\over 2}}(r) \over r^{1\over2} I_{k+{1\over2}}(R_s) },
\end{equation}
\begin{equation}\label{4.16}
  b_{kl}(r)= -c_{kl}\big[c_2-c_1(k^2+k)\big]
  \Big[{R_s^{1\over2} I_{k+{1\over 2}}(r)
   \over r^{1\over2} I_{k+{1\over2}}(R_s) }-{r^k\over R_s^k}\Big].
\end{equation}
  By (\ref{4.1}), (\ref{4.14}) and (\ref{4.16}) we compute
\begin{equation}\label{4.17}
\begin{array}{rll}
 & &\mathcal B(0)(b_{kl}(r)Y_{k,l}(\omega))+c_3c_{kl}Y_{k,l}(\omega)
  \\ [0.3 cm]
  & = &\displaystyle\big[b'_{kl}(R_s)+c_3c_{kl}\big]Y_{k,l}(\omega)
  \\ [0.3 cm]
  & = &\displaystyle \Big\{-\big[c_2-c_1(k^2+k)\big]{ I_{k+{3\over 2}}(R_s)
   \over  I_{k+{1\over2}}(R_s) }+c_3\Big\}c_{kl}Y_{k,l}(\omega)
  \\[0.3 cm]
  & = &\displaystyle  \Big \{-\mu\Big[{\bar\sigma\gamma\over R_s^2}
  (1-{1\over2}(k^2+k))-{1\over3}\tilde\sigma R_s\Big]{ I_{k+{3\over 2}}(R_s)
   \over  I_{k+{1\over2}}(R_s) }-\mu\bar\sigma(1-{\gamma\over R_s})
  +\mu\tilde\sigma\Big\}c_{kl}Y_{k,l}(\omega)
  \\ [0.5 cm]
  & =: \,& \displaystyle \Lambda_k(\gamma)c_{kl}Y_{k,l}(\omega).
\end{array}
\end{equation}
   Note that
$ I_{1/2}(r)=\sqrt{2\over \pi r}\sinh r $.
    One can easily verify that (\ref{2.3}) is equivalent to
\begin{equation}\label{4.18}
 f(R_s)=(1-{\gamma\over R_s}){I_{3/2}(R_s)
  \over R_s I_{1/2}(R_s)}\qquad
  \mbox{and}\qquad
  \tilde \sigma=3\bar\sigma(1-{\gamma\over R_s})
  {I_{3/2}(R_s) \over R_s I_{1/2}(R_s)}.
\end{equation}
  By using (\ref{4.18}), we can rewrite
\begin{equation}\label{4.19}
  \Lambda_k(\gamma)={\mu \bar\sigma\over R_s}
  \Big[ \gamma (h_k+j_k)- j_kR_s\Big]\qquad
  \mbox{for}\;\;\gamma>0,\;\;k=0,1,2\cdots,
\end{equation}
  where
\begin{equation}\label{4.20}
   h_{k}:=({k^2+k\over2}-1){I_{k+{3\over2}}(R_s)\over
   R_sI_{k+{1\over2}}( R_s)},
\qquad
  j_k:=1
  - {3I_{3\over2}(R_s)\over R_s I_{1\over 2}( R_s)}
  -{I_{3\over 2}(R_s)I_{k+{3\over 2}}(R_s)\over
   I_{1\over 2}( R_s)I_{k+{1\over 2}}(R_s)}.
\end{equation}
  By (4.46) in [\ref{wu-cui-07}] we have
\begin{equation}\label{4.21}
   j_0<0,\qquad j_1=0 \qquad\mbox{and}\qquad j_k>0\quad\mbox{for}\;\;  k\ge2.
\end{equation}
  It follows that
\begin{equation}\label{4.22}
  \Lambda_1(\gamma)\equiv 0
  \qquad \mbox{for}\;\;\gamma>0.
\end{equation}
  A direct computation also shows that
$$
  f'(R_s)=-{1\over R_s^2}\big[\gamma (h_0+j_0)-j_0 R_s\big],
$$
  so that by (\ref{2.4}) and (\ref{4.19}) we have
\begin{equation}\label{4.23}
  \Lambda_0(\gamma)= -\mu \bar\sigma R_s f'(R_s)\neq 0
  \qquad \mbox{for}\;\;\gamma>0.
\end{equation}

  From (\ref{4.4}) and the deduction (\ref{4.17}), we have the following result:

\medskip

  {\bf Lemma 4.1} \ \  {\em For any $\eta\in
  C^\infty(\Bbb S^2)$ with expansion $\displaystyle\eta=\sum_{k=0}^
  \infty\sum_{l=-k}^{k}c_{kl}Y_{k,l}(\omega)$, there holds
\begin{equation}\label{4.24}
  \partial\Psi(0)\eta=\sum_{k=0}^\infty
  \sum_{l=-k}^{k}\Lambda_k(\gamma)c_{kl}Y_{k,l}(\omega),
\end{equation}
  where $\Lambda_k(\gamma)$ is given by $(\ref{4.19})$
  with $\Lambda_0(\gamma) = -\mu\bar\sigma R_s f'(R_s)$
  and  $\Lambda_1(\gamma) \equiv 0$.
}

\medskip

\medskip

   Lemma 4.1 implies that for  $k\in \{0,1,2,\cdots\}$ and $\gamma>0$,
   $\Lambda_k(\gamma)$ is an eigenvalue of the operator $\partial\Psi(0)$.
   As mentioned before, the spectrum of $\partial\Psi(0)$, denoted by
   $\sigma(\partial\Psi(0))$,  consists of
   all eigenvalues, hence we have the following

\medskip

{\bf Corollary 4.2} \ \ {\em The spectrum of $\partial\Psi(0)$ is given by
$$
  \sigma(\partial\Psi(0))=\{\Lambda_k(\gamma);k=0,1,2,\cdots\}.
$$

}

\medskip

  Next, we study the property of eigenvalue $\Lambda_k(\gamma)$.
  We introduce
\begin{equation}\label{4.25}
  \gamma_k:={j_k\over h_k+j_k}R_s \qquad
  \mbox{for}\;\; k=2, 3, \cdots,
\end{equation}
  where $h_k$ and $j_k$ are given by (\ref{4.20}).

\medskip

  {\bf Lemma 4.3} \ \ {\em   For each $k\in\{2,3,\cdots\}$, we have  $\gamma_k>0$.
   Moreover,  $\displaystyle\lim_{k\to +\infty}\gamma_k=0$.

  }
\medskip

  {\bf Proof}. \ \  By (\ref{4.20}) and (\ref{4.21}),
$$
  h_k>0\qquad\mbox{and}\qquad 0<j_k<1\qquad\mbox{for}\;\;k=2,3,\cdots.
$$
  Thus we have
$$
   0<\gamma_k={j_k R_s\over h_k+j_k}<{R_s\over h_k+1}
   \qquad\mbox{for}\;\;k=2,3,\cdots.
$$
   Recall the well-known formula (cf. [\ref{wat-44}])
$$
  I_m(r)=\sqrt {1\over2\pi m}({e r\over 2m})^m
  \left( 1+O({1\over m})\right)
  \qquad \hbox{as}\;\; m\to +\infty.
$$
  We have
$$
  {I_{k+{3\over2}}(r)\over I_{k+{1\over2}}(r)}= e r
  {(2k+1)^{k+1}\over(2k+3)^{k+2}}
  \Big(1+O\Big(\frac{1}{k}\Big)\Big)=\frac{r}{2k}
  +O\Big(\frac{1}{k^2}\Big) \qquad \mbox{as}\;\; k\to +\infty,
$$
  and
$$
  {R_s\over h_k +1 }=\Big[({k^2+k\over2}-1){I_{k+{3\over2}}(R_s)\over
   R_sI_{k+{1\over2}}( R_s)}+1\Big]^{-1}R_s={4\over k}R_s+O({1\over k^2})
   \qquad \mbox{as}\;\; k\to +\infty.
$$
  It immediately follows that  $\displaystyle \lim_{k\to+\infty}\gamma_k=0$.
 \qquad $\Box$

\medskip

  Let
\begin{equation}\label{4.26}
  \gamma_*:=\sup\{\gamma_k;  k=2,3,\cdots\}.
\end{equation}
  By Lemma 4.3, we see that $\gamma_*>0$.
\medskip

 Now, we state the following useful properties of eigenvalues of $\partial\Psi(0)$ :

\medskip

{\bf Lemma 4.4} \ \ {\em $(i)$ For any $\gamma>0$, we have
  $\displaystyle\lim_{k\to+\infty}\Lambda_k(\gamma)=+\infty$.

  $(ii)$ If $\gamma>\gamma_*$, then $\Lambda_k(\gamma)>0$
  for all $k\in\{2,3,\cdots\}$; Moreover, 0 is an eigenvalue of
  $\partial\Psi(0)$ with geometric multiplicity $3$.

  $(iii)$ If  $0<\gamma<\gamma_*$, then there exists at least
  one integer $k_0\ge2 $ such that $\Lambda_{k_0}(\gamma)<0$.
  }

\medskip

  {\bf Proof}. \ \
  By (\ref{4.19}) and (\ref{4.25}), we see that $\Lambda_k(\gamma)$ can be
  expressed as
\begin{equation}\label{4.27}
   \Lambda_k(\gamma)={\mu\bar\sigma\over R_s}(h_k+j_k)(\gamma - \gamma_k)
    \qquad\mbox{for}\;\;k=2, 3,\cdots.
\end{equation}
  Recall that $h_k+j_k>0$ for all $k=2,3,\cdots$. By the proof of Lemma 4.3, we
  can show that
$$
  h_k+j_k={k\over 4}+O(1)\qquad \mbox{as}\;\;k\to+\infty.
$$
  By $\displaystyle\lim_{k\to+\infty}\gamma_k=0$, we immediately
  get the assertion $(i)$.

  Clearly,  the assertions $(ii)$ and $(iii)$ follow from
  (\ref{4.26}), (\ref{4.27}) and Lemma 4.1.  \qquad $\Box$

\medskip

{\bf Remark 4.5} \ \ $0\in \sigma(\partial\Psi(0))$ is due to the translation
  invariance of original free boundary problem (\ref{1.1}).  In fact, recall that
  $(\sigma_s^{[x_0]},p_s^{[x_0]},\rho_s^{[x_0]})$ is the translated radial stationary
  solution induce by $[x\to x+x_0]$ as introduced before. Note that
  $\rho_s^{[0]}=0$. By the reduction in section 2, we have $\Psi(\rho_s^{[x_0]})=0$
  for $|x_0|$ is sufficiently small. Let ${\bf r}_x:={x/ |x|}$ for $x\neq0$.
  By (2.3) in [\ref{fri-rei-02}],  we have ${\bf r}_x=(e_1, e_2, e_3)$ with
  $e_1, e_2, e_3\in \mbox{span}\{Y_{1,1}, Y_{1,0}, Y_{1,-1}\}$.
   Since $|x|=R_s+\rho_s^{[x_0]}$ is equivalent to $|x-x_0|=R_s$,
   we have
$$
  \rho_s^{[x_0]}=|x_0|\langle {\bf r}_{x_0}, {\bf r}_x\rangle + O(|x_0|^2)
  \qquad\mbox{as }\;\; |x_0|\to 0,
$$
  and $\langle {\bf r}_{x_0}, {\bf r}_x\rangle\in \mbox{span}\{Y_{1,1}, Y_{1,0}, Y_{1,-1}\}$.
  It implies that $\partial\Psi(0)Y_{1,l}=0$, $l=-1,0,1$. Moreover,
  for given $\delta>0$ sufficiently small,
$$
  \mathfrak M:=\{\rho_s^{[x_0]}, x_0\in \Bbb R^3, |x_0|<\delta\}
$$
  is a 3-dimensional Banach manifold with the tangent space at $0$ given
  by  $\mbox{span}\{Y_{1,1}, Y_{1,0}, Y_{1,-1}\}$.

\medskip
\hskip 1em

\section{Asymptotic stability}
\setcounter{equation}{0}
\hskip 1em

  In this section we study asymptotic behavior of transient solutions
  of problem (\ref{2.15}) and give the  proof of Theorem 1.2. Since
  $0\in\sigma(\partial\Psi(0))$,  the standard linearized stability
  theorem for parabolic differential equations in Banach spaces
  is not applicable. We have to employ a so-called {\em generalized
  principle of linearized stability} developed in [\ref{pru-sim-09}]
  to overcome this difficulty.

\medskip

 {\bf Theorem 5.1}\ \ {\em $(i)$ If $f'(R_s)<0$ and $\gamma>\gamma_*$,
  then the equilibrium  $0$ of problem $(\ref{2.15})$ is asymptotically
  stable in the following sense: There exists constant $\epsilon>0$
  such that for any $\rho_0\in \mathcal O_\delta$ satisfying
  $\|\rho_0\|_{h^{4+\alpha}(\Bbb S^2)}<\epsilon$, problem
  $(\ref{2.15})$ has a unique global classical solution $\rho(t)$ in
  $\mathcal O_\delta$, which converges exponentially fast
  to a translated equilibrium $\rho_s^{[x_0]}\in \mathfrak M$
  in $h^{4+\alpha}(\Bbb S^2)$ for some $x_0\in \Bbb R^3$,
  as $t\to +\infty$.

  $(ii)$ If $f'(R_s)>0$ or $0<\gamma<\gamma_*$,
  the equilibrium  $0$ is unstable.
}

\medskip

{\bf Proof}. \ \ $(i)$ Set $X_0:= h^{3+\alpha}(\Bbb S^2)$ and
  $X_1:=h^{4+\alpha}(\Bbb S^2)$. For given $T>0$, let $J:=[0,T)$
  and
$$
\begin{array}{l}
  \Bbb E_0(J):=BUC(J,X_0),
\\ [0.3 cm]
  \Bbb E_1(J):=BUC^1(J,X_0)\cap BUC(J,X_1).
\end{array}
$$
  Note that $X_1$ is densely embedded into $X_0$ and the
  little H\"{o}lder spaces have the following well-known
  interpolation property
\begin{equation}\label{5.1}
  (h^{\sigma_0}(\Bbb S^2),h^{\sigma_1}(\Bbb S^2))_{\theta,\infty}^0=
  h^{(1-\theta)\sigma_0+\theta\sigma_1}(\Bbb S^2) \qquad
  \mbox{if}\;\;(1-\theta)\sigma_0+\theta\sigma_1\notin \Bbb N,
\end{equation}
  where $\theta\in (0,1)$, $0<\sigma_0<\sigma_1$ and
  $(\cdot,\cdot)_{\theta,\infty}^0$
  denotes the continuous interpolation functor (cf. [\ref{amann},
  \ref{lunardi}]). We see that $X_1$ is the trace space of $\Bbb E_1(J)$.

  Set ${\mathscr A}:= \partial\Psi(0)$. By Lemma 3.2,  ${\mathscr A}
   \in \mathcal H(X_1,X_0)$. Then by Remark III. 3.4.2 (b) in
  [\ref{amann}] and (\ref{5.1}), we have
\begin{equation}\label{5.2}
  (\Bbb E_0(J),\Bbb E_1(J)) \mbox{ is a pair of maximal
  continuous regularity for }  \mathscr A.
\end{equation}

  For $f'(R_s)<0$ and $\gamma>\gamma_*$, by (\ref{4.23})
  and Lemma 4.4, we have $ \Lambda_0(\gamma)=-\mu\bar\sigma R_sf'(R_s)>0$,
  $\Lambda_k(\gamma)>0$, $k\ge2$, and $0$ is an eigenvalue of geometric
  multiplicity 3.  By Lemma 4.4 $(i)$, we see that there exists constant
  $\varpi>0$ such that
\begin{equation}\label{5.3}
   \sigma(\mathscr A)\backslash \{0\}\subset \Bbb C+:=
   \{z\in \Bbb C:  {\rm Re} z>\varpi\}.
\end{equation}
  Denote the kernel  and the range of $\mathscr A$ by
  $N(\mathscr A)$ and $R(\mathscr A)$, respectively.
  Then by Lemma 4.1, for $\gamma>\gamma_*$ we have
\begin{equation}\label{5.4}
  N(\mathscr A)=\mbox{span}\{Y_{1,1}, Y_{1,0}, Y_{1,-1}\}
  \qquad\mbox{and}\qquad
 N(\mathscr A)\oplus R(\mathscr A)=X_0.
\end{equation}

  By Remark 4.5,  we see that $\mathfrak M=\{\rho_s^{[x_0]},
  x_0\in \Bbb R^3, |x_0|<\delta\}$ is a 3-dimensional
  $C^1$-manifold near the equilibrium $0$ and the tangent
  space for $\mathfrak M$ at $0$ is given by $N(\mathscr A)$.
  Thus from this observation and together with (\ref{5.2})--(\ref{5.4}),
  we see that all assumptions of Theorem 3.1 in [\ref{pru-sim-09}]
  are satisfied, so that the desired assertion
  $(i)$ follows from this theorem immediately.   More precisely,
  there exists $\epsilon>0$ such that the solution $\rho(t)$ of
  problem (\ref{2.15}) exists globally for any $\rho_0\in
  \mathcal O_\delta$ satisfying
  $\|\rho_0\|_{h^{4+\alpha}(\Bbb S^2)}<\epsilon$,  and
  moreover,  there exist constant  $M>0$ and a point
  $x_0\in\Bbb R^3$ close to $0$, such that
\begin{equation}\label{5.5}
   \|\rho(t)-\rho_s^{[x_0]}\|
  _{h^{4+\alpha}(\Bbb S^2)}
  \le  M e^{-\varpi t}\|\rho_0\|_{h^{4+\alpha}(\Bbb S^2)}\qquad
  \mbox{for}\;\;
  t\ge0.
\end{equation}

  $(ii)$ If $f'(R_s)>0$ or $0<\gamma<\gamma_*$, then by
  (\ref{4.23}) and Lemma 4.4 $(iii)$,
  we see that $\sigma(-\mathscr A)\cap\{z\in\Bbb C: {\rm Re}
  z>0\}$ is not empty. It follows from Theorem 9.1.3 in
  [\ref{lunardi}]  that the equilibrium $0$ is unstable.
  The proof is complete. \qquad $\Box$

\medskip

  {\bf Remark 5.2}\ \ Theorem 5.1 $(i)$ can be also proved by using
  a center manifold theorem for quasi-invariant parabolic differential
  equations in Banach spaces developed in [\ref{cui-09}].
  Since original free boundary problem (\ref{1.1}) is invariant
  under coordinate translations,  the reduce differential equation
  (\ref{2.15}) possesses quasi-invariance under a natural local Lie
  group action induced by coordinate translations. With the aid of
  Theorem 2.1 in [\ref{cui-09}], we can show that $\mathfrak M$
  is a center manifold which attracts nearby transient solutions
  and get the desired assertion. The proof is more complicated
  though it can obtain more delicate information than the proof given
  here,  we refer interested readers to the proof of Theorem 1.1
  in [\ref{cui-09}], see also the proof of main result in [\ref{wu-cui-09}]
  and [\ref{wu-zhou-13}] for similar free boundary problems.

\medskip

  Now, we give the proof of our main result Theorem 1.2.

\medskip

  {\bf Proof of Theorem 1.2.} \ \  Recall that for $0<\tilde\sigma/\bar\sigma<\theta_*$,
  the free boundary problem (\ref{1.1}) has two radial stationary solution
  $(\sigma_s,p_s,\rho_s)$ with radius $R_{s1}$ and $R_{s2}$,
  for $0<R_{s1}<R_{s2}$. By (\ref{2.4}), we have $f'(R_{s1})>0$ and
  $f'(R_{s2})<0$. Thus by Lemma 2.1, Lemma 2.2 and Theorem 5.1 $(i)$,
  we see that for the radial stationary solution $(\sigma_s,p_s,\rho_s)$
  with radius $R_{s2}$ in case $\gamma>\gamma_*$,
  there exists constant $\epsilon>0$ such that for any
  $\rho_0\in \mathcal O_\delta$  satisfying
  $\|\rho_0\|_{h^{4+\alpha}(\Bbb S^2)}<\epsilon$,
  free boundary problem $(\ref{1.1})$ has a unique
  global solution $(\sigma(t),p(t), \rho(t))$ which is given by
$$
  \sigma(t)=\Theta_*^{\rho(t)}\mathcal U(\rho(t)),\qquad
  p(t)=\Theta_*^{\rho(t)}\mathcal V(\rho(t)),
$$
  and $\rho(t)$ is the solution of problem (\ref{2.15}) with
  $\rho(0)=\rho_0$, where $\mathcal U$ and $\mathcal V$
  are given by (\ref{2.12}). By (\ref{5.5}) and the reduction in
  section 2, we see that there exists $x_0\in \Bbb R^3$ such
  that  $(\sigma(t),p(t), \rho(t))$ converges exponentially fast
  to  $(\sigma_s^{[x_0]}, p_s^{[x_0]}, \rho_s^{[x_0]})$ in
  $h^{2+\alpha}(\Omega_{\rho(t)})\times h^{4+\alpha}
  (\Omega_{\rho(t)}) \times h^{4+\alpha}(\Bbb S^2)$,
  as time goes to infinity.

  Similarly by Lemma 2.1, Lemma 2.2 and
  Theorem 5.1 $(ii)$, we see that
  the radial stationary solution $(\sigma_s,p_s,\rho_s)$
  with radius $R_{s2}$ is unstable for $0<\gamma<\gamma_*$,
  and the radial stationary
  solution $(\sigma_s,p_s,\rho_s)$ with radius $R_{s1}$
  is unstable for all $\gamma>0$.  This completes the proof.
  \qquad $\Box$

\medskip
\hskip 1em

  {\bf Acknowledgement.} \ \  This work is supported by the
  National Natural Science Foundation of China under grant 11371147
  and PAPD of Jiangsu Higher Education Institutions.

\vsss

\end{document}